\documentclass[11pt]{article}
\usepackage{amsfonts,amssymb,amsmath,latexsym}
\usepackage[english]{babel}
\usepackage[latin1]{inputenc}
\newcommand{\Nz}{\mathbb{N}}
\newcommand{\Rz}{\mathbb{R}}
\newcommand{\Cz}{\mathbb{C}}
\newcommand{\Kz}{\mathbb{K}}
\newcommand{\Hz}{\mathbb{H}}
\newcommand{\isomorph}{\overset{\sim}=}
\newcommand{\trace}{\mathsf{tr\,}}
\newcommand{\qR}{\mathcal R}

\newcommand{\com}[2]{\left[#1,#2 \right]}
\newenvironment {PROOF}{\textsc{Proof:} \small}{{\hspace*{\fill} $\square$}}
\newtheorem {LEMMA} {Lemma} [section]

\newtheorem {THEOREM} [LEMMA] {Theorem}

\newtheorem {DEFINITION}[LEMMA] {Definition}

\title {Solutions of the Quantum-Yang-Baxter-Equation \\
from Symmetric Spaces}

\author{Martin Bordemann\thanks{M.Bordemann@univ-mulhouse.fr}\\
          Laboratoire de Math\'{e}matiques et Applications\\
          Facult\'{e} des Sciences et Techniques\\
          Universit\'{e} de Haute Alsace, Mulhouse\\
          4, rue des Fr\`{e}res Lumi\`{e}re\\
          68093 Mulhouse, France \\[3ex]
          Martin Walter\thanks{Martin.Walter@physik.uni-freiburg.de}\\
          Fakult\"{a}t f\"{u}r Physik\\Universit\"{a}t Freiburg \\
          Hermann-Herder-Str. 3 \\
          79104 Freiburg i.~Br., Germany\\[3ex]
          FR-TH 01/06\\}

\date{July 2001}

\begin {document}

\maketitle
\thispagestyle{empty}

\begin {abstract}
 We show that for each semi-Riemannian locally symmetric space the curvature
 tensor gives rise to a rational solution $r$ of the classical Yang-Baxter
 equation with spectral parameter. For several Riemannian globally
 symmetric spaces $M$ such as real, complex and quaternionic Grassmann
 manifolds we explicitly compute solutions $R$ of the quantum Yang Baxter
 equations (represented in the tangent spaces of $M$) which generalize
 the quantum $R$-matrix found by Zamolodchikov and Zamolodchikov in 1979.
\end {abstract}

\section{Introduction and Results}\label{IntroSec}
Since their discoveries not only the classical
Yang-Baxter equation (CYBE) but also the
the quantum Yang-Baxter equation (QYBE)
are playing an important r\^ole in several
branches of physics and mathematics.
Following Sklyanin, the CYBE can be formulated for an element
$r\in\mathfrak{g}\wedge \mathfrak{g}$ on
an abstract Lie algebra $\mathfrak{g}$ and takes
the form
\[
[r_{12},r_{13}] + [r_{12},r_{23}]+ [r_{13},r_{23}] = 0~.
\]
We use the standard notation from the theory of
Hopf algebras, where subscripts describe on which
places in an $n$-fold tensor product the quantities
in question act. We refer to solutions of this equation
as classical $r$-matrices. On an associative algebra
${\cal A}$ we can formulate the quantum version
of above equation for an element $\qR\in{\cal A}\otimes
{\cal A}$ as
\begin{equation}
\qR_{12}\qR_{13}\qR_{23} = \qR_{23}\qR_{13}\qR_{12}~.
\end{equation}
Solutions of this equation are commonly called
quantum $\qR$-matrices. For a precise treatment
see the important work of Drinfel'd \cite{drinfeld} or
the standard literature on quantum groups
(e.g. \cite{chari} and references therein).
As it occurs in  statistical models for example,
these equations possibly depend on additional
spectral parameters.

A special case of
interest  are quantum-$\qR$-matrices which
can be written as a deformation of classical
$r$-matrices:
\begin{equation}\label{clqr}
\qR = \mathsf{id} + \hbar r + o(\hbar^2)~.
\end{equation}
What we are actually looking for are quantum
$\qR$-matrices fulfilling this relation
for some given classical $r$-matrices.

The classical $r$-matrices we have in mind
are arising from locally symmetric spaces. It
can be shown that curvature tensor of these
spaces obey CYBE with spectral parameter
(theorem \ref{motiv}). This is because
there exists a Casimir element of a reductive
Lie algebra $\mathfrak{k}$ of the isotropy group
$K$ of a symmetric space obeying CYBE.
The curvature tensor is exactly this Casimir element
being represented in the tangent space of the symmetric space
(theorem \ref{motivreason}). The problem is now
how to obtain
a deformation in the sense of \eqref{clqr}
fulfilling QYBE for a representation
of given $\mathfrak{k}$ in the tangent space. It turns out that
a heuristic  ansatz in quantities one
can associate to a symmetric space
in a canonical manner yield such solutions
(section 3). Using this we explicitly
calculate for some of the symmetric spaces
which has been classified in \cite{helgason}
quantum $\qR$-matrices. In particular
for spaces of constant sectional curvature,
projective spaces and Grassmann manifolds
the ansatz turned out to be successful
(see section 4).

\section{Symmetric Spaces and CYBE}
Let us briefly recall the geometric background
necessary for our considerations (see e.g. \cite{helgason}
and \cite{sakai} for more details):

 A semi-Riemannian
differentiable manifold $(M,g)$
 is called {\em locally symmetric} iff around
every point $m$ the geodesic symmetry ($x\mapsto -x$
in an exponential chart) is a local isometry. An equivalent
criterion is the statement that the curvature tensor
$R$ of its Levi-Civita connection $\nabla$ is
covariantly constant, $\nabla R=0$. Standard examples
are the $n$-spheres, but also all Riemann surfaces
equipped with the standard constant curvature metric.
A semi-Riemannian locally symmetric space is called
{\em globally symmetric} or {\em symmetric} iff the geodesic
symmetry around every point extends to a global isometry
of $M$. This is for instance no longer the case for Riemann
surfaces of genus $\geq 2$. Semi-Riemannian globally symmetric
spaces are known to be homogeneous under its Lie group of all
isometries whence they admit a group theoretic representation
$M=G/K$ where $G$ is a Lie group acting isometrically on $M$ and $K$ is
the closed subgroup of all those transformations of $G$ fixing
a chosen point $p$. Consider a pair $(G,K)$ of a Lie group $G$ and
a closed subgroup
$K$ with its associated pair of Lie algebras $(\mathfrak{g},\mathfrak{k})$.
$(G,K)$ is called a {\em semi-Riemannian symmetric pair} iff i) there is an
involutive automorphism $\sigma$ of $G$ such that $K$ is a subgroup
of the group $G_{\sigma}$ of all the fixed points of $\sigma$ and
contains the identity component of $G_\sigma$, and ii) there is a
nondegenerate symmetric bilinear form $(~,~)$ on $\mathfrak{m}$,
the eigenspace associated to the eigenvalue $-1$ of
$\theta:=T_e\sigma$, which is invariant under $\theta$ and the restriction
of the adjoint representation of $G_\sigma$ to $\mathfrak{m}$. Note
that $\mathfrak{m}$ is identified with the tangent space $T_pM$
and that the Lie subalgebra $\mathfrak{k}$ is the eigenspace associated
to the eigenvalue $+1$ of $\theta$ and that $\theta$ is an involutive
automorphism of the Lie algebra $\mathfrak{g}$.
There are  the well-known commutation
relations:

\begin{equation}\label{stvertrel}
\com{\mathfrak{k}}{\mathfrak{k}} \subset \mathfrak{k}~,~~~
\com{\mathfrak{k}}{\mathfrak{m}} \subset \mathfrak{m}~,~~~
\com{\mathfrak{m}}{\mathfrak{m}} \subset \mathfrak{k}~.
\end{equation}
Such a pair $(\mathfrak{g},\mathfrak{k})$ will be also
called symmetric pair for short. If $M=G/K$ is a Riemannian
globally symmetric space then the pair can be chosen in such a
way that $K$ be a compact subgroup of $G$.

The most known examples of symmetric spaces are
spaces of constant sectional curvature (like
the sphere $S^n$ and the hyperbolic space $H^n$),
the complex and quaternionic projective spaces
$\Cz$P$^n$ and $\Hz$P$^n$, and the Grassmann manifolds
$SO(p+q)/S(O(p)\times O(q))$, $SU(p+q)/S(U(p)\times U(q))$
and $Sp(p+q)/Sp(p)\times Sp(q)$.

We are going to prove in two ways that the curvature tensor of a symmetric
space fulfills CYBE. The first one is based
on direct differential geometric calculations, the second will make this
fact obvious for group theoretical reasons.

Let $(M,g)$ a Riemannian manifold with metric tensor $g$
and Levi-Civita connection $\nabla$. For the second
covariant derivative (along two vector fields $X,Y$ on $M$)
the equation \cite{ko+no}
\begin{equation}
\nabla^2_{XY} = \nabla_X \nabla_Y - \nabla_{\nabla_X Y}~.
\end{equation}
is valid. Therefore we obtain for
any $(r,s)$-tensor field  $T\in {\cal T}_s^r(M)$ on $M$
that its antisymmetric part of the second covariant
derivative obeys
\begin{equation}\label{crucial R1}
(\nabla^2_{XY} -\nabla^2_{YX}) T = R(X,Y)\circ T~,
\end{equation}
or in local coordinates (with
${R_{klm}}^i\frac{\partial}{\partial x^i}:=R(\frac{\partial}{\partial
x^k},\frac{\partial}{\partial x^l})\frac{\partial}{\partial x^m}$):
\begin{equation}\label{crucial R2}
T_{j_1\dots j_s;lk}^{i_1\dots i_r} - T_{j_1\dots j_s;kl}^{i_1\dots i_r} =
\sum_\alpha R_{klm}^{\hspace{3ex}i_\alpha}
T_{j_1.........j_s}^{i_1\dots m \dots i_r}
-\sum_\alpha R_{klj_\alpha}^{\hspace{3ex}m}
T^{i_1.........i_r}_{j_1\dots m\dots j_s}~.
\end{equation}
We are now able to formulate the following theorem.
\begin{THEOREM}\label{motiv}
Let $M$ be a semi-Riemannian locally symmetric space,  $R_{ijkl}$ the components
of the curvature tensor. Then, at each point $p\in M$,
\begin{equation}\label{classical r matrix}
\hat{r}:= \frac{1}{\lambda_1-\lambda_2}
R^{i\;\;k}_{\;\;j\;\;l}\, E_i^{\;j}
\otimes E_k^{\;l}~\in \mathsf{End}(T_pM)\otimes \mathsf{End}(T_pM)
\otimes \Cz(\lambda_1,
\lambda_2)
\end{equation}
is a rational solution of the CYBE with spectral parameters.
\end{THEOREM}
\begin{PROOF}
The fact that $R$ is covariant constant is essential
for the proof. Inserting $R$ in equation \eqref{crucial R2}
and rearranging the indices by using the metric, one finds
\begin{equation}\label{ltsgleichung}
R^{i \;\;k }_{\;\;m \;\;l } R^{m \;\;r }_{\;\;n \;\;s } -
R^{m \;\;k }_{\;\;n \;\;l } R^{i \;\;r }_{\;\;m \;\;s } +
R^{i \;\;k }_{\;\;n \;\;m } R^{m \;\;r }_{\;\;l \;\;s }-
R^{ i\;\;m }_{\;\;n \;\;l } R^{k \;\;r }_{\;\;m \;\; s} =0~.
\end{equation}
This is exactly the condition for the vanishing of the commutators
one has to calculate in CYBE with \eqref{classical r matrix}
as classical $r$-matrix.
\end{PROOF}\\

This key observation was the starting point of our
considerations. For semi Riemannian globally symmetric
spaces this theorem can also be proven by group theoretic methods
using representation theory of the Casimir element.

Let $M$ be a globally symmetric space,
$(\mathfrak{g},\mathfrak{k})$ be the associated
symmetric pair to $M\isomorph G/K$ and $\rho$ the adjoint
representation of $\mathfrak{k}$ on $\mathfrak{m}$.
Suppose $\mathfrak{k}$ is a reductive Lie algebra,
i.e. there exists a nondegenerate, symmetric and
invariant bilinear form on $\mathfrak{k}$ for which
we can define the corresponding Casimir element $t$.
Identifying $\mathfrak{k}$ with its dual via the
bilinear form, and introducing an orthonormal
basis $\{e_\alpha\}$, we can write $t=\sum_{\alpha}e_\alpha\otimes
e_\alpha$. This element gives rise to rational
solutions of CYBE with spectral parameters (cf. \cite{drinfeld}).
Under these assumptions we can show
the following theorem which
was also proven in \cite{wang}.
\begin{THEOREM}\label{motivreason}
The curvature tensor $R$ of a locally symmetric space is the
in $\mathfrak{m}\isomorph T_pM$ represented Casimir element $t$,
\begin{equation}
R(p) = (\rho\otimes \rho)(t)~.
\end{equation}
\end{THEOREM}
\begin{PROOF}
The curvature tensor of a symmetric space at a point $p\in M$
can be written as
\[
(R(X,Y)Z)(p) = -[[X,Y],Z](p)
\]
($\forall X,Y,Z\in \mathfrak{m}$). Then the statement follows
from the fact that there exists an $\mathsf{Ad} K$-invariant,
non-degenerate scalar product on $\mathfrak{g}$ and can be easily
worked out using an orthonormal basis.
\end{PROOF}\\

From this point of view theorem \ref{motiv} becomes obvious
since the solutions of CYBE obtained there are just because
of the special representation we have chosen.

\section{Quantum-$\qR$-Matrices associated to Symmetric Pairs}
In this section we will shortly explain the main
ideas how to obtain solutions of the QYBE from the
curvature tensor of a symmetric space. We assume to have
the following ingredients, where $\mathfrak{g}$ denotes a
finite dimensional Lie algebra over the field $\Kz$:
\begin{enumerate}
\item On $\mathfrak{g}$ there exists an invariant,
symmetric, nondegenerate bilinear form $\kappa$.
\item There exists an involutive automorphism $\theta:\mathfrak{g}\to
\mathfrak{g}$ splitting the Lie Algebra into a direct sum
of vector spaces $\mathfrak{g}=\mathfrak{k}\oplus \mathfrak{m}$
with commutations relations \eqref{stvertrel}.
\item We have a faithful representation $\rho:\mathfrak{g}\to {\cal A}$
in some associative algebra ${\cal A}$.
\end{enumerate}
Because of $i.)$ we can identify $\mathfrak{g}$ with its dual
and get the isomorphism $(\forall X,Y,Z\in\mathfrak{g})$
\begin{eqnarray*}
\Phi:\mathfrak{g}\otimes\mathfrak{g} & \to & \mathsf{End}(\mathfrak{g})~,\\
\Phi(X\otimes Y)(Z) & = & X\kappa(Y,Z)~.
\end{eqnarray*}
Using $ii.)$ and $iii.)$ we define the following
objects.
\begin{DEFINITION}
Under above assumptions define the elements:
\begin{equation}
\begin{array}{lccl}
C  & := & \Phi^{-1}(\mathsf{id}) & \in   \mathfrak{g}\otimes\mathfrak{g}~, \\
G  & := & \Phi^{-1}(\theta)  & \in  \mathfrak{g}\otimes\mathfrak{g}~,
\end{array}
\end{equation}
and their pendants represented in ${\cal A}$:
\begin{equation}
\begin{array}{lccl}
\hat{C}  & := & (\rho\otimes\rho)C  & \in  {\cal A}\otimes{\cal A}~, \\
\hat{G}  & := & (\rho\otimes\rho)G  & \in  {\cal A}\otimes{\cal A}~.
\end{array}
\end{equation}
\end{DEFINITION}
Since it is easily verified that $C$ is the so-called
quadratic Casimir element of $\mathfrak{g}$
and is therefore $\mathfrak{g}$-invariant,
it obeys the CYBE with spectral parameters (see e.g. \cite{drinfeld}).
Although the sum of $C$ and $G$ is not $\mathfrak{g}$-
but $\mathfrak{k}$-invariant it yields a solution
of this equation (cf. e.g. \cite{BFLS}, where all solutions
of CYBE with spectral parameters depending on $C$ and $G$
have been classified.) by setting
\begin{equation}\label{classr}
r:= \frac{C+G}{\lambda_1-\lambda_2}
~\in~\mathfrak{k}\otimes\mathfrak{k}\otimes \Cz(\lambda_1,\lambda_2)~.
\end{equation}
The classical $r$-matrix defined above has universal character
because it is defined in the abstract Lie algebra $\mathfrak{k}$ and
does not depend on any representation. For that
reason it is clear that we obtain matrix solutions
$\hat{r}=(\rho\otimes\rho) r \in {\cal A}\otimes{\cal A}\otimes
\Cz(\lambda_1,\lambda_2)$ for any representation
$\rho$.

In order to maintain the aspect that the classical
$r$-matrix corresponds the curvature tensor we
choose $\mathfrak{m}$ as representation space for $\mathfrak{g}$,
i.e. ${\cal A}=\mathsf{End}(\mathfrak{m})$. We are actually looking
for $\qR$-matrices obeying the QYBE, which
can be written as a deformation of the classical $r$-matrix $\hat{r}$:
\begin{equation}\label{qrgeneral}
\hat{\qR} = \mathsf{id} + \hbar \hat{r}+
\sum_{k=2}^\infty \hbar^k\hat{\qR}^{(k)}~,\qquad
\hat{\qR}^{(k)}
~\in~\mathsf{End}(\mathfrak{m})
\otimes \mathsf{End}(\mathfrak{m})\otimes \Cz(\lambda_1,\lambda_2)~,
\end{equation}
where $\hbar$ is an arbitrary parameter. The aim is
to find explicit expressions for the quantities
$\hat{\qR}^{(k)}$, where we make --guided by the classical
solution-- a heuristic ansatz, that the solutions
\eqref{qrgeneral} of QYBE can be written as
\begin{equation}\label{qrprecise}
\hat{\qR} = \mathsf{id} + a(\lambda_1-\lambda_2,\hbar)\hat{G}
                        + b(\lambda_1-\lambda_2,\hbar)\hat{C}~,
\end{equation}
where $a,b\in\Cz(\lambda_1,\lambda_2)[[\hbar]]$ have to be
determined. One should keep in mind that this ansatz is
by no means unique, but in some cases successful as we show below.

However, in general it is not possible
to realize the whole Lie algebra $\mathfrak{g}$ on
the vector space $\mathfrak{m}$, whereas the Lie subalgebra
$\mathfrak{k}$ is canonical represented by the adjoint representation.
This rather technical problem can be avoided if we introduce
a new symmetric pair $(\tilde{\mathfrak{g}},\mathfrak{k})$ with
an involutive automorphism $\tilde{\theta}$ such that
there exists a representation of $\tilde{\mathfrak{g}}$ on
$\mathfrak{m}$. Therefore, one has to replace $\theta$ by $\tilde{\theta}$
in the above definition and the resulting $\hat{C}$ and $\hat{G}$
constitute the $\qR$-matrix for which we have to find the
coefficients $a$ and $b$.

\section{Examples}
We will present examples (cf. \cite{walter}) for which
we could successfully determine the coefficients
for solutions of \eqref{qrprecise}. Instead of proving
everything in detail we just sketch the proof. Taking
equation \eqref{qrprecise} and inserting it in QYBE,
we obtain restrictions for the coefficients, which
can be solved. The other way around, one easily but
lengthly checks by direct calculations that the given
$\qR$-matrices fulfill the QYBE.

In the succeeding examples we use the following notation.
By $E_{ij}$ we denote the elementary matrices of
given dimension having a $1$ in the $i$-th row and $j$-th
column and being $0$ elsewhere. All vector spaces and their
tensor products are considered over the reals.  For $\rho$ we choose
the usual matrix multiplication.

\begin{enumerate}
\item Spaces of constant curvature. $k$ takes the
value $1,0,-1$ depending whether it is a space
of positive, zero or negative sectional curvature.
\begin{flalign*}
&\begin{array}{l@{\quad=\quad}l}
\tilde{\mathfrak{g}}  &  \mathfrak{gl}(n,\Rz)\\
\mathfrak{k}          &  \mathfrak{so}(n)\\
\mathfrak{m}          &  \Rz^n\\
\tilde{\theta}(X)        &  -X^T\\
\kappa(X,Y)           &  \trace(\rho(X)\rho(Y))\\
\hat{G}               &  -E_{ij}\otimes E_{ij}\\
\hat{C}               &  E_{ij}\otimes E_{ji}\\
a                     &  \frac{\hbar}{\lambda_1-\lambda_2
                         +\frac{k(n-2)}{2}\hbar}\\
b                     &  \frac{\hbar}{\lambda_1-\lambda_2}
\end{array}&
\end{flalign*}

In order to compute the coefficients $a$ and $b$ we made use
of the following identities satisfied by $\hat{C}$ and $\hat{G}$
which are readily verified:

\begin{enumerate}
\item ${\hat{G}}_{12}{\hat{G}}_{23} = - {\hat{C}}_{13}{\hat{G}}_{23}
  = - {\hat{G}}_{12}{\hat{C}}_{13} = - {\hat{G}}_{12}{\hat{G}}_{13}{\hat{G}}_{23}
        = 1/n {\hat{G}}_{12}{\hat{C}}_{13}{\hat{G}}_{23}
  = - {\hat{G}}_{12}{\hat{C}}_{13}{\hat{C}}_{23} $ \\$
        = - {\hat{C}}_{23}{\hat{G}}_{13}{\hat{C}}_{12}
  = - {\hat{C}}_{12}{\hat{C}}_{13}{\hat{G}}_{23}
        $~,
\item ${\hat{G}}_{23}{\hat{G}}_{12} = - {\hat{G}}_{23}{\hat{C}}_{13}
  = - {\hat{C}}_{13}{\hat{G}}_{12} = - {\hat{G}}_{23}{\hat{G}}_{13}{\hat{G}}_{12}
        = 1/n {\hat{G}}_{23}{\hat{C}}_{13}{\hat{G}}_{12}
  = -  {\hat{G}}_{23}{\hat{C}}_{13}{\hat{C}}_{12} $ \\$
        = - {\hat{C}}_{12}{\hat{G}}_{13}{\hat{C}}_{23}
  = - {\hat{C}}_{23}{\hat{C}}_{13}{\hat{G}}_{12}
        $~,
\item ${\hat{C}}_{12}{\hat{G}}_{13}{\hat{G}}_{23}
         = {\hat{G}}_{23}{\hat{G}}_{13}{\hat{C}}_{12} = -{\hat{G}}_{23}$~,
\item ${\hat{G}}_{12}{\hat{G}}_{13}{\hat{C}}_{23}
         = {\hat{C}}_{23}{\hat{G}}_{13}{\hat{G}}_{12} = -{\hat{G}}_{12}$~,
\item ${\hat{C}}_{12}{\hat{C}}_{13}{\hat{C}}_{23}
         = {\hat{C}}_{23}{\hat{C}}_{13}{\hat{C}}_{12} = {\hat{C}}_{13}$~,
\item ${\hat{G}}_{13}{\hat{G}}_{23}
         = -{\hat{C}}_{12}{\hat{G}}_{23} =-{\hat{G}}_{13}{\hat{C}}_{12}$~,
\item ${\hat{G}}_{12}{\hat{G}}_{13}
         = -{\hat{C}}_{23}{\hat{G}}_{13}= -{\hat{G}}_{12}{\hat{C}}_{23}$~,
\item ${\hat{G}}_{13}{\hat{G}}_{12}
         = -{\hat{C}}_{23}{\hat{G}}_{12} = -{\hat{G}}_{13}{\hat{C}}_{23}$~,
\item ${\hat{G}}_{23}{\hat{G}}_{13}
         = -{\hat{C}}_{12}{\hat{G}}_{13} = -{\hat{G}}_{23}{\hat{C}}_{12}$~,
\item ${\hat{C}}_{12}{\hat{C}}_{13}
         = {\hat{C}}_{23}{\hat{C}}_{12} = {\hat{C}}_{13}{\hat{C}}_{23}$~,
\item ${\hat{C}}_{13}{\hat{C}}_{12}
         = {\hat{C}}_{12}{\hat{C}}_{23} = {\hat{C}}_{23}{\hat{C}}_{13}$~.
\end{enumerate}

The solution for the sphere $S^n$, i.e. $k=1$, is already presented in
\cite{reshetikhin} and was originally computed by \cite{zandz}.

\item The complex projective space $\Cz$P$^n$.
\begin{flalign*}
&\begin{array}{l@{\quad=\quad}l}
\tilde{\mathfrak{g}}  &  \mathfrak{gl}(n,\Cz)\\
\mathfrak{k}          &  \mathfrak{u}(n)\\
\mathfrak{m}          &  \Cz^n\\
\tilde{\theta} (X)       &  -X^\dag\\
\kappa(X,Y)           &  {\cal R}e\,\trace(\rho(X)\rho(Y))\\
\hat{G}               &  -E_{ij}\otimes E_{ij}- iE_{ij}\otimes iE_{ij} \\
\hat{C}               &  E_{ij}\otimes E_{ji}- iE_{ij}\otimes iE_{ji} \\
a                     &  \frac{\hbar}{\lambda_1-\lambda_2
                         +n\hbar}\\
b                     &  \frac{\hbar}{\lambda_1-\lambda_2}
\end{array}&
\end{flalign*}

Here the identities (c)--(k) of the preceding example remain wheras (a) and
(b) are replaced by:
\begin{enumerate}
\item ${\hat{G}}_{12}{\hat{G}}_{23} = - {\hat{C}}_{13}{\hat{G}}_{23}
  = - {\hat{G}}_{12}{\hat{C}}_{13}
        = 1/(2n) {\hat{G}}_{12}{\hat{C}}_{13}{\hat{G}}_{23}
        $~,
\item ${\hat{G}}_{23}{\hat{G}}_{12} = - {\hat{G}}_{23}{\hat{C}}_{13}
  = - {\hat{C}}_{13}{\hat{G}}_{12}
        = 1/(2n) {\hat{G}}_{23}{\hat{C}}_{13}{\hat{G}}_{12}
        $~,
\end{enumerate}

\item The quaternionic projective space $\Hz$P$^n$.
\begin{flalign*}
&\begin{array}{l@{\quad=\quad}l}
\tilde{\mathfrak{g}}  &  \mathfrak{gl}(n,\Hz)\\
\mathfrak{k}          &  \mathfrak{sp}(n)\\
\mathfrak{m}          &  \Hz^n\\
\tilde{\theta}(X)        &  -\bar{X}^T\\
\kappa(X,Y)           &  {\cal S}c\,\trace(\rho(X)\rho(Y))\\
\hat{G}               &  -E_{ij}\otimes E_{ij}- iE_{ij}\otimes iE_{ij}
                         -jE_{ij}\otimes jE_{ij}- kE_{ij}\otimes kE_{ij}\\
\hat{C}               &  E_{ij}\otimes E_{ji}- iE_{ij}\otimes iE_{ji}
                         -jE_{ij}\otimes jE_{ji}- kE_{ij}\otimes kE_{ji} \\
a                     &  \frac{\hbar}{\lambda_1-\lambda_2
                         +(2n+2)\hbar}\\
b                     &  \frac{\hbar}{\lambda_1-\lambda_2}
\end{array}&
\end{flalign*}

Here again the identities (c)--(k) of the preceding example remain wheras (a) and
(b) are replaced by:
\begin{enumerate}
\item ${\hat{G}}_{12}{\hat{G}}_{23} = - {\hat{C}}_{13}{\hat{G}}_{23}
  = - {\hat{G}}_{12}{\hat{C}}_{13} = (1/2){\hat{G}}_{12}{\hat{G}}_{13}{\hat{G}}_{23}
        = 1/(4n) {\hat{G}}_{12}{\hat{C}}_{13}{\hat{G}}_{23}
  = (1/2) {\hat{G}}_{12}{\hat{C}}_{13}{\hat{C}}_{23} $ \\$
        = (1/2) {\hat{C}}_{23}{\hat{G}}_{13}{\hat{C}}_{12}
  = (1/2) {\hat{C}}_{12}{\hat{C}}_{13}{\hat{G}}_{23}
        $~,
\item ${\hat{G}}_{23}{\hat{G}}_{12} = - {\hat{G}}_{23}{\hat{C}}_{13}
  = - {\hat{C}}_{13}{\hat{G}}_{12} = (1/2){\hat{G}}_{23}{\hat{G}}_{13}{\hat{G}}_{12}
        = 1/(4n) {\hat{G}}_{23}{\hat{C}}_{13}{\hat{G}}_{12}
  = (1/2)  {\hat{G}}_{23}{\hat{C}}_{13}{\hat{C}}_{12} $ \\$
        = (1/2) {\hat{C}}_{12}{\hat{G}}_{13}{\hat{C}}_{23}
  = (1/2) {\hat{C}}_{23}{\hat{C}}_{13}{\hat{G}}_{12}
        $~,
\end{enumerate}
\end{enumerate}

After having introduced solutions for projective spaces
the following three examples we encounter unfortunately possess
a less geometric meaning. Let $p,q\in\Nz$ and set
$P=\{1,\dots, p\}$, $Q=\{p+1,\dots,p+q\}$. We do not indicate the
identities for $ \hat{C}$ and $\hat{G}$ (as in the preceding examples);
they are computed similarly.
\begin{enumerate}
\item[$iv.)$] The symmetric pair
$(\mathfrak{gl}(p+q,\Rz),\mathfrak{gl}(p)\times\mathfrak{gl}(q))$.
Defining
\[
I_{p,q} = \begin{pmatrix} -\mathbf{1}_p & 0 \\ 0 &  \mathbf{1}_q
\end{pmatrix}~,
\]
where $\mathbf{1}_p$ respectively  $\mathbf{1}_q$ is the $p$-
respectively $q$-dimensional identity matrix, we find
\begin{flalign*}
&\begin{array}{l@{\quad=\quad}l}
\tilde{\mathfrak{g}}  &  \mathfrak{gl}(p+q,\Rz) \\
\mathfrak{k}          &  \mathfrak{gl}(p)\times\mathfrak{gl}(q) \\
\mathfrak{m}          &  \Rz^{(p+q)}\\
\tilde{\theta} (X)       &  I_{p,q}XI_{p,q}\\
\kappa(X,Y)           &  \trace(\rho(X)\rho(Y))\\
\hat{G}               &  \sum_{i,j\in P}+\sum_{i,j\in Q}-
                         \sum_{\substack{i\in P\\j\in Q }}
                         -\sum_{\substack{i\in Q \\j\in P}}
                          E_{ij}\otimes E_{ji}\\
\hat{C}               &   \sum_{i,j\in P\cup Q} E_{ij}\otimes E_{ji} \\
a                     &  \frac{\hbar}{\lambda_1-\lambda_2}\\
b                     &  \frac{\hbar}{\lambda_1-\lambda_2}
\end{array}&
\end{flalign*}
\item[$v.)$]  The symmetric pair
$(\mathfrak{gl}(p+q,\Rz),\mathfrak{so}(p,q))$.
\begin{flalign*}
&\begin{array}{l@{\quad=\quad}l}
\tilde{\mathfrak{g}}  &  \mathfrak{gl}(p+q) \\
\mathfrak{k}          &  \mathfrak{so}(p,q) \\
\mathfrak{m}          &  \Rz^{(p+q)}\\
\tilde{\theta}(X)        &  -I_{p,q}X^TI_{p,q}\\
\kappa(X,Y)           &  \trace(\rho(X)\rho(Y))\\
\hat{G}               &  - \sum_{i,j\in P} - \sum_{i,j\in Q}+
                         \sum_{\substack{i\in P\\j\in Q }}
                         + \sum_{\substack{i\in Q\\j\in P}}
                          E_{ij}\otimes E_{ij} \\
\hat{C}               &  \sum_{i,j\in P\cup Q} E_{ij}\otimes E_{ji} \\
a                     &  \frac{\hbar}{\lambda_1-\lambda_2+
                         \frac{p+q-2}{2}\hbar} \\
b                     &  \frac{\hbar}{\lambda_1-\lambda_2}
\end{array}&
\end{flalign*}
\item[$vi.)$] The symmetric pair
$(\mathfrak{gl}(2n,\Rz),\mathfrak{gl}(n,\Cz))$.
Introducing
\[
I  = \begin{pmatrix} 0 & \mathbf{1} \\ -\mathbf{1} & 0
\end{pmatrix}~,
\]
we get the following solution of QYBE.
\begin{flalign*}
&\begin{array}{l@{\quad=\quad}l}
\tilde{\mathfrak{g}}  &  \mathfrak{gl}(2n,\Rz) \\
\mathfrak{k}          &  \mathfrak{gl}(n,\Cz) \\
\mathfrak{m}          &  \Rz^{2n} \\
\tilde{\theta}(X)     &  IXI \\
\kappa(X,Y)           &  \trace(\rho(X)\rho(Y))\\
\hat{G}               &  \sum_{i,j}^n E_{i+n,j+n} \otimes E_{j,i}
                         + E_{i,j} \otimes E_{j+n,i+n}
                         - E_{i,j+n} \otimes E_{j,i+n}
                         - E_{i+n,j} \otimes E_{j+n,i} \\
\hat{C}               &   \sum_{i,j=1}^{2n} E_{ij}\otimes E_{ji} \\
a                     &  \frac{\hbar}{\lambda_1-\lambda_2}\\
b                     &  \frac{\hbar}{\lambda_1-\lambda_2}
\end{array}&
\end{flalign*}
\end{enumerate}
All examples considered so far have in common that they make
use of the trick mentioned in the last section. The
$\qR$-matrices are unitary up to a factor as readily verified
and could been determined by the ansatz \eqref{qrprecise}.
Moreover, expanding the coefficients $a$ into geometric
series and comparing it with \eqref{qrgeneral} we can
can read off explicitly the $\hat{\qR}^{(k)}$.

The next examples we would like to discuss are
Grassmann manifolds. Here our ansatz fails.
Nevertheless we immediately see, that the
solutions we obtained for projective spaces yield solutions
of QYBE for these spaces. We will concentrate
on the real Grassmann manifolds of the form $SO(p+q)/S(O(p)\times O(q))$
since the complex and quaternionic cases can be
treated in a similar way. The corresponding Lie algebra
$\mathfrak{g}=\mathfrak{so}(p+q)$ splits into
\begin{eqnarray*}
\mathfrak{k} & = &
\left\{ \begin{pmatrix} A & 0 \\ 0 & D \end{pmatrix}
\bigg|\, A\in \mathfrak{so}(p)~, D\in\mathfrak{so}(q) \right\}
 \isomorph\mathfrak{so}(p)\times \mathfrak{so}(q)~,\\
\mathfrak{m} & = &
\left\{ \begin{pmatrix} 0  & B \\ -B^T & 0 \end{pmatrix}
\bigg| \,B \mbox{ arbitrary $(p\times q)$-matrix} \right\}\isomorph
\mathsf{End}(\Rz^q,\Rz^p)~.
\end{eqnarray*}
The canonical representation of $\mathfrak{k}$  on $\mathfrak{m}$
is given by
\begin{eqnarray}\label{grassmannrep}
\rho: \mathfrak{k} & \to & \mathsf{End}(\mathfrak{m})~,\\
\rho(A,D)(B) &:= & AB-BD \nonumber
\end{eqnarray}
Let us recall the solution we found in our first example $i.)$
for different dimensions $p$ and $q$. Namely, we get
\begin{eqnarray*}
\qR^p & = &  \mathsf{id}+\frac{\hbar}{\lambda_1-\lambda_2+\frac{p-2}{2}\hbar}
\hat{G}_p +
\frac{\hbar}{\lambda_1-\lambda_2}\hat{C}_p~, \\
\qR^q & = &  \mathsf{id}+\frac{\hbar}{\lambda_1-\lambda_2+\frac{q-2}{2}\hbar}
\hat{G}_q +
\frac{\hbar}{\lambda_1-\lambda_2}\hat{C}_q~.
\end{eqnarray*}
We immediately can formulate
\begin{THEOREM}
Let $M=SO(p+q)/S(O(p)\times O(q))$ and $\hat{t}$ the
in $\mathfrak{m}\isomorph T_pM$ represented Casimir
element. Let $v\in\mathfrak{m}\otimes\mathfrak{m}$
be an arbitrary element. Then
\begin{equation}
\hat{\qR}(v):= \qR^p v \qR^q
\end{equation}
is a solution of the QYBE whose leading terms in the
power series expansion in the parameter $\hbar$ take
the form
\begin{equation}
\hat{\qR} = \mathsf{id}+ \frac{\hbar}{\lambda-\mu}\hat{t}+
\frac{\hbar^2}{(\lambda-\mu)^2}\left(\frac{1}{2}\hat{t}^2-\mathsf{id}\right)
+o(\hbar^3)~.
\end{equation}
\end{THEOREM}
\begin{PROOF}
The first assumption is obvious, since
$\qR^p$ and $\qR^q$ are solutions of QYBE themselves. We
then obtain
\begin{eqnarray*}
\hat{\qR}_{12}\hat{\qR}_{13}\hat{\qR}_{23}(v)
& = & \qR^p_{12}\qR^p_{13}\qR^p_{23}\, v\,\qR^q_{23}\qR^q_{13}\qR^q_{12}\\[1ex]
& = & \qR^p_{23}\qR^p_{13}\qR^p_{12}\, v\,\qR^q_{12}\qR^q_{13}\qR^q_{23}\\[1ex]
& = & \hat{\qR}_{23}\hat{\qR}_{13}\hat{\qR}_{12}(v)~.
\end{eqnarray*}
The second part can be proven by  straightforward calculation.
\end{PROOF}\\

In the same manner we find the solutions for the
remaining Grassmann manifolds there the classical
$r$-matrix corresponds to the curvature tensor.

\section{Acknowledgments}
We would like to thank W. Soergel for making us
aware of the references \cite{reshetikhin,wang}.

\end{document}